\documentclass[11pt]{amsart}
\usepackage{amssymb,latexsym}
\usepackage{pstcol,pstricks,color}

 \numberwithin{equation}{section}

\newtheorem{theorem}{Theorem}[section]
\newtheorem{proposition}[theorem]{Proposition}
\newtheorem{corollary}[theorem]{Corollary}

\newtheorem{remark}[theorem]{Remark}

\newtheorem{lemma}[theorem]{Lemma}

\newtheorem{example}[theorem]{Example}

\newcommand{\dis}{\displaystyle}

\title{Bell Polynomials and $k$-generalized Dyck Paths}

\begin{document}
\maketitle
\begin{center}
Toufik Mansour$^\dag$ and Yidong Sun\footnote{Corresponding author:
Yidong Sun, sydmath@yahoo.com.cn.}$^\ddag$

$^\dag$Department of Mathematics, University of Haifa, 31905 Haifa,
Israel\\
$^\dag$Center for Combinatorics, LPMC, Nankai University, 300071
Tianjin, P.R. China

$^\ddag$Department of Mathematics, Dalian Maritime University, 116026 Dalian, P.R. China\\[5pt]

{\it toufik@math.haifa.ac.il, sydmath@yahoo.com.cn}
\end{center}\vskip0.5cm

\subsection*{Abstract}
A {\em k-generalized Dyck path} of length $n$ is a lattice path from
$(0,0)$ to $(n,0)$ in the plane integer lattice
$\mathbb{Z}\times\mathbb{Z}$ consisting of horizontal-steps $(k, 0)$
for a given integer $k\geq 0$, up-steps $(1,1)$, and down-steps
$(1,-1)$, which never passes below the $x$-axis. The present paper
studies three kinds of statistics on $k$-generalized Dyck paths:
"number of $u$-segments", "number of internal $u$-segments" and
"number of $(u,h)$-segments". The Lagrange inversion formula is used
to represent the generating function for the number of
$k$-generalized Dyck paths according to the statistics as a sum of
the partial Bell polynomials or the potential polynomials. Many
important special cases are considered leading to several surprising
observations. Moreover, enumeration results related to $u$-segments
and $(u,h)$-segments are also established, which produce many new
combinatorial identities, and specially, two new expressions for
Catalan numbers.

\medskip

{\bf Keywords}: Bell polynomials, Potential polynomials, $k$-paths,
Catalan numbers

\noindent {\sc 2000 Mathematics Subject Classification}: Primary
05A05, 05A15; Secondary 05C90

\section{Introduction}

Let $\mathfrak{L}_{n,k}$ denote the set of lattice paths of length
$n$ from $(0,0)$ to $(n,0)$ in the plane integer lattice
$\mathbb{Z}\times\mathbb{Z}$ consisting of horizontal-steps $h=(k,
0)$ for a given integer $k\geq 0$, up-steps $u=(1,1)$, and
down-steps $d=(1,-1)$. Let $\mathfrak{L}_{n,k}^{m,j}$ be the set of
lattice paths in $\mathfrak{L}_{n,k}$ with $m$ up-steps and $j$
horizontal-steps. Let $L$ be any lattice path in
$\mathfrak{L}_{n,k}^{m,j}$. A {\em $u$-segment} of $L$ is a maximal
sequence of consecutive up-steps in $L$. Define $\alpha_i(L)$ to be
the number of $u$-segments of length $i$ in $L$ and call $L$ having
the $u$-segments of type $1^{\alpha_1(L)}2^{\alpha_2(L)}\cdots$. Let
$\mathfrak{L}_{n,k,r}^{m,j}$ be the subset of lattice paths in
$\mathfrak{L}_{n,k}^{m,j}$ with $r$ $u$-segments.

A {\em k-generalized Dyck path} or {\em k-path} (for short) of
length $n$ is a lattice path in $\mathfrak{L}_{n,k}$ which never
passes below the $x$-axis. By our notation, a {\em Dyck path} is a
0-path, a {\em Motzkin path} is a 1-path and a {\em Schr\"oder path}
is a 2-path. Let $\mathfrak{P}_{n,k}^{m,j}$ denote the set of
$k$-paths of length $n$ (i.e. $n=2m+kj$) with $m$ up-steps and $j$
horizontal-steps and let $\mathfrak{Q}_{n,k}^{m,j}$ be the subset of
$k$-paths in $\mathfrak{P}_{n,k}^{m,j}$ with no horizontal-step at
$x$-axis. Define $\mathfrak{P}_{n,k,r}^{m,j}\
(\mathfrak{Q}_{n,k,r}^{m,j})$ to be the subset of $k$-paths in
$\mathfrak{P}_{n,k}^{m,j}\ (\mathfrak{Q}_{n,k}^{m,j})$ with $r$
$u$-segments.

In \cite{mansun}, we study two kinds of statistics on Dyck paths:
"number of $u$-segments" and "number of internal $u$-segments". In
this paper, we consider these two statistics together with "number
of $(u,h)$-segments" in the more extensive setting of $k$-paths. In
order to do this, we present two necessary tools : Lagrange
inversion formula and the potential polynomials.

\subsection*{Lagrange Inversion Formula~\cite{wilf}}
If $f(x)=\sum_{n\geq 1}f_nx^n$ with $f_1\neq 0$, then the
coefficients of the composition inverse $g(x)$ of $f(x)$ $({\rm
namely}, f(g(x))=g(f(x))=x)$ can be given by
\begin{eqnarray}\label{eqn 1.1}
[x^n]g(x)=\frac{1}{n}[x^{n-1}]\big(\frac{x}{f(x)}\big)^{n}.
\end{eqnarray}
More generally, for any formal power series $\Phi(x)$,
\begin{eqnarray}\label{eqn 1.2}
[x^n]\Phi(g(x))=\frac{1}{n}[x^{n-1}]\Phi'(x)\big(\frac{x}{f(x)}\big)^{n},
\end{eqnarray}
for all $n\geq 1$, where $\Phi'(x)$ is the derivative of $\Phi(x)$
with respect to $x$.
\subsection*{The Potential Polynomials~\cite[pp. 141,157]{comtet}} The potential polynomials
$\textbf{P}_{n}^{(\lambda)}$  are defined for each complex number
$\lambda$ by
\begin{eqnarray*}
1+\sum_{n\geq
1}\textbf{P}_{n}^{(\lambda)}\frac{x^n}{n!}&=&\Big\{1+\sum_{n\geq
1}f_n\frac{x^n}{n!}\Big\}^\lambda,
\end{eqnarray*}
which can be represented by Bell polynomials
\begin{eqnarray}\label{eqn 1.3}
\textbf{P}_{n}^{(\lambda)}=\textbf{P}_{n}^{(\lambda)}(f_1,f_2,f_3,\dots)=\sum_{1\leq
k\leq n}\binom{\lambda}{k}k!{\bf{B}}_{n,k}(f_1,f_2,f_3,\dots),
\end{eqnarray}
where  ${\bf B}_{n,i}\big(x_1,x_2,\cdots\big)$ is the partial Bell
polynomial on the variables $\{x_j\}_{j\geq 1}$ (see~\cite{Be}).

In this paper, with the Lagrange inversion formula, we can represent
the generating functions for the number of $k$-paths according to
our statistics (see Sections 2-4) as a sum of partial Bell
polynomials or the potential polynomials. For example,
\begin{eqnarray*}
\sum_{P\in\mathfrak{P}_{n,k}^{m,j}}\prod_{i\geq1}t_i^{\alpha_i(P)}
&=&\frac{1}{m+1}{\binom{m+j}{m}}\sum_{r=0}^{m}{\binom{m+j+1}{r}}\frac{r!}{m!}{\bf
B}_{m,r}\big(1!t_1,2!t_2,\cdots\big), \\
\sum_{Q\in
\mathfrak{Q}_{n,k}^{m,j}}\prod_{i\geq1}t_i^{\alpha_i(Q)}&=&
\frac{m}{(m+j+1)(m+j)}\binom{m+j}{j}\frac{{\bf{P}}_m^{(m+j+1)}(1!t_1,2!t_2,\cdots)}{m!}.
\end{eqnarray*}
We consider a number of important special cases. These lead to
several surprising results. Moreover, enumeration results related to
$u$-segments and $(u,h)$-segments are also established in Section 5,
producing many new combinatorial identities and in particular the
following two new expressions for the Catalan numbers:
\begin{eqnarray*}
\sum_{p=0}^{[n/2]}\frac{1}{2p+1}\binom{3p}{p}\binom{n+p}{3p}&=&\frac{1}{n+1}\binom{2n}{n},\hskip0.5cm (n\geq 0)\\
\sum_{p=0}^{[(n-1)/2]}\frac{1}{2p+1}\binom{3p+1}{p+1}\binom{n+p}{3p+1}&=&\frac{1}{n+1}\binom{2n}{n},\hskip0.5cm
(n\geq 1).
\end{eqnarray*}

\section{"$u$-segments" statistics in $k$-paths}
We start this section by studying the generating function for the
number of $k$-paths of length $n$ according to the statistics
$\alpha_1,\alpha_2,\ldots$, that is,
\begin{eqnarray*}
P(x,z;{\bf t})=P(x,z;t_1,t_2,\ldots)=\sum_{m,j\geq0}x^mz^j
\sum_{P\in\mathfrak{P}_{n,k}^{m,j}}\prod_{i\geq1}t_i^{\alpha_i(P)}.
\end{eqnarray*}

\begin{proposition}\label{pro2.1}
The ordinary generating function $P(x,z;{\bf t})$ is given by
\begin{eqnarray}\label{eqn 2.1}
P(x,z;{\bf t})=1+zP(x,z;{\bf t})+\sum_{j\geq1}t_jx^{j}P^j(x,z;{\bf
t})+z\sum_{j\geq1}t_jx^{j}P^{j+1}(x,z;{\bf t}).
\end{eqnarray}
\end{proposition}
\begin{proof}
Note that $P(x,z;{\bf t})$ can be written as $P(x,z;{\bf
t})=1+zP(x,z;{\bf t})+\sum_{j\geq1}P_j(x,z;{\bf t})$, where
$P_j(x,z;{\bf t})$ is the generating function for the number of
$k$-paths with initial $u$-segment of length $j$ according to the
statistics $\alpha_1,\alpha_2,\ldots$. An equation for $P_j(x,z;{\bf
t})$ is obtained from the first return decomposition of a $k$-path
starting with a $u$-segment of length $j$: either $P=u^jdP^{(1)}d
P^{(2)} d\ldots P^{(j-1)}d P^{(j)}\mbox{ or }P= u^jhP^{(1)} d
P^{(2)} d\ldots P^{j}dP^{(j+1)},$ where $P^{(1)},\ldots,P^{(j)}$ are
$k$-paths, see Figure~\ref{fDD2.1}.

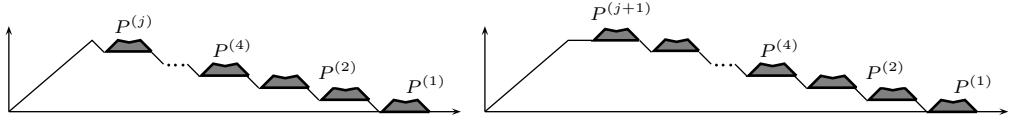
\begin{figure}[h]
\begin{pspicture}(6.2,2)
\psline[unit=18pt,linewidth=.5pt]{->}(0,0)(9.5,0)\psline[unit=18pt,linewidth=.5pt]{->}(0,0)(0,1.8)
\psline[unit=18pt,linewidth=.5pt](0,0)(1.75,1.5)(2,1.25)
\put(1.15,.65){\psline[unit=17pt,linewidth=1pt,fillstyle=solid,fillcolor=gray](1.25,.25)(.25,.25)(.5,.5)(.75,.45)(1,.5)(1.25,.25)}
\psline[unit=18pt,linewidth=.5pt](3,1.25)(3.25,1)\pscircle*[unit=18pt,linewidth=.5pt](3.36,.98){0.03}
\pscircle*[unit=18pt,linewidth=.5pt](3.5,.98){0.03}\pscircle*[unit=18pt,linewidth=.5pt](3.66,.98){0.03}
\psline[unit=18pt,linewidth=.5pt](3.75,1)(4,.75)
\multiput(-1.38,-.18)(.8,-.16){4}{\put(3.78,.5){\psline[unit=18pt,linewidth=1pt,fillstyle=solid,fillcolor=gray](1.25,.25)(.25,.25)(.5,.5)(.75,.45)(1,.5)(1.25,.25)}}
\multiput(0,0)(.8,-.16){3}{\psline[unit=18pt,linewidth=.5pt](5,.75)(5.25,.5)}
\put(1.4,1.05){\tiny$P^{(j)}$} \put(2.7,.73){\tiny$P^{(4)}$}
\setlength\unitlength{10pt}
\rput*[r](4.65,.55){\tiny$P^{(2)}$}\put(15,.7){\tiny$P^{(1)}$}
\end{pspicture}
\begin{pspicture}(6.2,2)
\psline[unit=18pt,linewidth=.5pt]{->}(0,0)(11,0)\psline[unit=18pt,linewidth=.5pt]{->}(0,0)(0,1.8)
\psline[unit=18pt,linewidth=.5pt](0,0)(1.75,1.5)(2.25,1.5)
\put(1.3,.8){\psline[unit=17pt,linewidth=1pt,fillstyle=solid,fillcolor=gray](1.25,.25)(.25,.25)(.5,.5)(.75,.45)(1,.5)(1.25,.25)}
\psline[unit=18pt,linewidth=.5pt](3.25,1.5)(3.5,1.25)\put(1.4,1.22){\tiny$P^{(j+1)}$}
\put(2.1,.66){\psline[unit=17pt,linewidth=1pt,fillstyle=solid,fillcolor=gray](1.25,.25)(.25,.25)(.5,.5)(.75,.45)(1,.5)(1.25,.25)}
\put(.95,0){\psline[unit=18pt,linewidth=.5pt](3,1.25)(3.25,1)\pscircle*[unit=18pt,linewidth=.5pt](3.36,.98){0.03}
\pscircle*[unit=18pt,linewidth=.5pt](3.5,.98){0.03}\pscircle*[unit=18pt,linewidth=.5pt](3.66,.98){0.03}
\psline[unit=18pt,linewidth=.5pt](3.75,1)(4,.75)
\multiput(-1.38,-.18)(.8,-.16){4}{\put(3.78,.5){\psline[unit=18pt,linewidth=1pt,fillstyle=solid,fillcolor=gray](1.25,.25)(.25,.25)(.5,.5)(.75,.45)(1,.5)(1.25,.25)}}
\multiput(0,0)(.8,-.16){3}{\psline[unit=18pt,linewidth=.5pt](5,.75)(5.25,.5)}
\put(2.7,.73){\tiny$P^{(4)}$} \setlength\unitlength{10pt}
\rput*[r](4.65,.55){\tiny$P^{(2)}$}\put(15,.7){\tiny$P^{(1)}$}}
\end{pspicture}

\caption{First return decomposition of a $k$-path starting with
exactly $j$ up-steps.}\label{fDD2.1}
\end{figure}

Thus $P_j(x,z;{\bf t})=t_{j}x^{j}P^{j}(x,z;{\bf
t})+zt_{j}x^{j}P^{j+1}(x,z;{\bf t})$. Hence, $P(x,z;{\bf t})$
satisfies $P(x,z;{\bf t})=1+zP(x,z;{\bf
t})+\sum_{j\geq1}t_jx^{j}P^j(x,z;{\bf
t})+z\sum_{j\geq1}t_jx^{j}P^{j+1}(x,z;{\bf t})$, as required.
\end{proof}

Define $y=y(x,z;{\bf t})=xP(x,z;{\bf t})$ and $T(x)=1+\sum_{j\geq
1}t_jx^j$. Then $(\ref{eqn 2.1})$ reduces to $y=(x+zy)T(y)$. Let
$y^*=y(x,zx;{\bf t})$, then we have $y^*=x(1+zy^*)T(y^*)$.

\begin{theorem}\label{theo 2.2}
For any integers $n,m\geq 1$ and $k,j\geq 0$,
\begin{eqnarray*}
\sum_{P\in\mathfrak{P}_{n,k}^{m,j}}\prod_{i\geq1}t_i^{\alpha_i(P)}
&=&\frac{1}{m+1}{\binom{m+j}{m}}\sum_{r=0}^{m}{\binom{m+j+1}{r}}\frac{r!}{m!}{\bf
B}_{m,r}\big(1!t_1,2!t_2,\cdots\big).
\end{eqnarray*}
\end{theorem}
\begin{proof}
Using (\ref{eqn 1.2}) and (\ref{eqn 1.3}), we obtain
\begin{eqnarray*}
\sum_{P\in\mathfrak{P}_{n,k}^{m,j}}\prod_{i\geq1}t_i^{\alpha_i(P)}&=&
[x^{m+j+1}z^j]y^*=\frac{[x^{m+j}z^j]}{m+j+1}\left\{(1+zx)T(x)\right\}^{m+j+1}\\
&=&
\frac{1}{m+j+1}{\binom{m+j+1}{j}}[x^{m}]T(x)^{m+j+1}\\
&=&\frac{1}{m+j+1}{\binom{m+j+1}{j}}\frac{{\bf
P}_{m}^{(m+j+1)}\big(1!t_1,2!t_2,\cdots\big)}{m!}\\
&=&\frac{1}{m+1}{\binom{m+j}{m}}\sum_{r=0}^{m}{\binom{m+j+1}{r}}\frac{r!}{m!}{\bf
B}_{m,r}\big(1!t_1,2!t_2,\cdots\big),
\end{eqnarray*}
as claimed.
\end{proof}

Replace $t_i$ by $qt_i$ in Theorem \ref{theo 2.2} and note that
\begin{eqnarray}
{\bf
B}_{m,r}\big(t_1,t_2,\cdots\big)&=&\sum_{\kappa_m(r)}\frac{m!}{r_1!r_2!\cdots
r_m!}\left(\frac{t_1}{1!}\right)^{r_1}\left(\frac{t_2}{2!}\right)^{r_2}\cdots
\left(\frac{t_m}{m!}\right)^{r_m},\label{eqn bella}\\
{\bf B}_{m,r}\big(1!qt_1,2!qt_2,\cdots\big)&=&q^r{\bf
B}_{m,r}\big(1!t_1,2!t_2,\cdots\big),\label{eqn bellb}
\end{eqnarray}
where the summation $\kappa_m(r)$ is for all the nonnegative integer
solutions of $r_1+r_2+\cdots+r_m=r$ and $r_1+2r_2+\cdots +mr_m=m$,
we have
\begin{corollary}\label{coro 2.3}
For any integers $n,m,r\geq 1$ and $k,j\geq 0$, there holds
\begin{eqnarray}\label{eqn 2.2}
\sum_{P\in\mathfrak{P}_{n,k,r}^{m,j}}\prod_{i\geq1}t_i^{\alpha_i(P)}
&=&\frac{1}{m+1}{\binom{m+j}{m}}{\binom{m+j+1}{r}}\frac{r!}{m!}{\bf
B}_{m,r}\big(1!t_1,2!t_2,\cdots\big).
\end{eqnarray}
\end{corollary}
By comparing the coefficient of $t_1^{r_1}t_2^{r_2}\cdots t_m^{r_m}$
in Corollary \ref{coro 2.3}, one can obtain that
\begin{corollary}
The number of $k$-paths in $\mathfrak{P}_{n,k,r}^{m,j}$ with
$u$-segments of type $1^{r_1}2^{r_2}\cdots m^{r_m}$ is
$\frac{1}{m+1}{\binom{m+j}{m}}{\binom{m+j+1}{r}}\binom{r}{r_1,r_2,\cdots,r_m}$.
Specially, the number of Dyck paths of length $2m$ with $u$-segments
of type $1^{r_1}2^{r_2}\cdots m^{r_m}$ is
$\frac{1}{m+1}{\binom{m+1}{r}}\binom{r}{r_1,r_2,\cdots,r_m}$, (the
case $k=0$ implies $j=0$).
\end{corollary}

\subsection{Applications}
In what follows we consider many special cases of $T(x)$. These
produce several interesting results, as described in
Examples~\ref{ex1}-\ref{ex11}. We also obtain several identities
involving Stirling numbers of the first (second) kind, idempotent
numbers and other combinatorial sequences.

\begin{example}\label{ex1} Let $T(x)=1+q(e^x-1)=(1-q)+qe^x$, that is, $t_i={q/i!}$ for all
$i\geq 1$. Note that ${\bf B}_{n,i}\big(q,q,q,\cdots\big)=S(n,i)q^i$
\cite[pp.135]{comtet}, where $S(n,i)$ is the Stirling number of the
second kind. Then Theorem \ref{theo 2.2} gives
\begin{eqnarray*}
\sum_{P\in\mathfrak{P}_{n,k}^{m,j}}\frac{m!\prod_{i\geq
1}q^{\alpha_i(P)}}{\prod_{i\geq1}(i!)^{\alpha_i(P)}}&=&
\frac{1}{m+j+1}\binom{m+j+1}{j}\sum_{i=0}^m\binom{m+j+1}{i}i!S(m,i)q^i\\
&=&\frac{1}{m+j+1}\binom{m+j+1}{j}\sum_{i=0}^m\binom{m+j+1}{i}i^mq^i(1-q)^{m+j-i+1},
\end{eqnarray*}
which leads to
$\sum_{P\in\mathfrak{P}_{n,k}^{m,j}}\frac{m!}{\prod_{i\geq1}(i!)^{\alpha_i(P)}}={\binom{m+j+1}{j}}(m+j+1)^{m-1}$
when $q=1$. By Corollary \ref{coro 2.3}, we have
\begin{eqnarray*}
\sum_{P\in\mathfrak{P}_{n,k,r}^{m,j}}\frac{m!}{\prod_{i\geq1}(i!)^{\alpha_i(P)}}&=&
\frac{1}{m+j+1}\binom{m+j+1}{j}\binom{m+j+1}{r}r!S(m,r).
\end{eqnarray*}
\end{example}

\begin{example} Let $T(x)=1+qxe^x$ which is equivalent to
$t_i={q/(i-1)!}$ for all $i\geq 1$. Note that ${\bf
B}_{m,i}\big(q,2q,3q,\cdots\big)=\binom{m}{i}i^{m-i}q^i$, which are
called the idempotent numbers \cite[pp.135]{comtet} when $q=1$. Then
Corollary \ref{coro 2.3} leads to
\begin{eqnarray*}
\dis\sum_{P\in\mathfrak{P}_{n,k,r}^{m,j}}\prod_{i\geq1}\left\{\frac{1}{(i-1)!}\right\}^{\alpha_i(P)}&=&
\frac{1}{m+j+1}{\binom{m+j+1}{j}}{\binom{m+j+1}{r}}\frac{r^{m-r}}{(m-r)!}.
\end{eqnarray*}
\end{example}

\begin{example} If $T(x)={(e^x-1)/x}$, then
$t_i={1/(i+1)!}$ for all $i\geq 1$. It is well known that the
Stirling numbers of the second kind satisfy
${\big(\frac{e^x-1}{x}\big)^k/k!}=\sum\limits_{m\geq
0}S(m+k,k){x^m/(m+k)!}$. Thus, Theorem \ref{theo 2.2} leads to
\begin{eqnarray*}
\dis\sum_{P\in\mathfrak{P}_{n,k}^{m,j}}\prod_{i\geq1}\frac{1}{\big((i+1)!\big)^{\alpha_i(P)}}
&=&\frac{(m+j)!}{(2m+j+1)!}\binom{m+j+1}{j}S(2m+j+1,m+j+1).
\end{eqnarray*}
\end{example}

\begin{example} If $T(x)=\frac{1}{x}\ln\frac{1}{1-x}$, then
$t_i={1/(i+1)}$ for all $i\geq 1$. It is well known that the
Stirling numbers of the first kind $s(n,i)$ satisfy
${\big(\frac{1}{x}\ln\frac{1}{1-x}\big)^k/k!}=\sum\limits_{m\geq
0}|s(m+k,k)|{x^m/(m+k)!}$. Thus, Theorem \ref{theo 2.2} leads to
\begin{eqnarray*}
\sum_{P\in\mathfrak{P}_{n,k}^{m,j}}\prod_{i\geq1}\frac{1}{(i+1)^{\alpha_i(P)}}&=&
\frac{(m+j)!}{(2m+j+1)!}\binom{m+j+1}{j}|s(2m+j+1,m+j+1)|.
\end{eqnarray*}
\end{example}

\begin{example}\label{ex 2.9} If $T(x)=1+q\ln\frac{1}{1-x}$, then
$t_i={q/i}$ for all $i\geq 1$. Using the fact that ${\bf
B}_{n,i}\big(0!q,1!q,2!q,\cdots\big)=|s(n,i)|q^i$
\cite[pp.135]{comtet} together with Corollary \ref{coro 2.3}, we
have
\begin{eqnarray*}
\sum_{P\in\mathfrak{P}_{n,k,r}^{m,j}}\prod_{i\geq1}\left\{\frac{1}{i}\right\}^{\alpha_i(P)}&=&
\frac{1}{m+1}{\binom{m+j}{m}}{\binom{m+j+1}{r}}\frac{r!}{m!}|s(m,r)|.
\end{eqnarray*}
\end{example}

\begin{example} If $T(x)={1/(1-x)^{\lambda}}$, then
$t_i=\binom{\lambda+i-1}{i}$ for all $i\geq 1$, where $\lambda$ is
an indeterminant. So, Theorem \ref{theo 2.2} leads to
\begin{eqnarray*}
\sum_{P\in\mathfrak{P}_{n,k}^{m,j}}{\prod_{i\geq1}\binom{\lambda+i-1}{i}^{\alpha_i(P)}}
&=&\frac{1}{m+1}\binom{m+j}{m}\binom{\lambda(m+j+1)+m-1}{m},
\end{eqnarray*}
which generates that when $\lambda=1$ the set
$\mathfrak{P}_{n,k}^{m,j}$ is counted by
$\frac{1}{m+1}\binom{m+j}{m}\binom{2m+j}{m}$, in particular,
$\mathfrak{P}_{n,k}^{m,m}$ is counted by
$\frac{1}{m+1}\binom{2m}{m}\binom{3m}{m}$.
\end{example}

\begin{example} Let $T(x)=1+x+x^2+\cdots+x^r$, that is,
$t_i=1$ for $1\leq i\leq r$ and $t_i=0$ for all $i\geq r+1$. Then
Theorem \ref{theo 2.2} gives
\begin{eqnarray*}
\sum_{P\in\mathfrak{P}_{n,k}^{m,j}}{\prod_{i\geq
1}^r1^{\alpha_i(P)}{\prod_{i\geq r+1}0^{\alpha_i(P)}}}
&=&\frac{1}{m+1}\binom{m+j}{m}\sum_{i=0}^{m+1}(-1)^i\binom{m+j+1}{i}\binom{2m+j-(r+1)i}{m+j},
\end{eqnarray*}
which implies that the number of $k$-paths $P$ of length $2n$ with
no $u$-segments of length greater than $r$ is given by
$$\frac{1}{m+1}\binom{m+j}{m}\sum_{i=0}^{m+1}(-1)^i\binom{m+j+1}{i}\binom{2m+j-(r+1)i}{m+j}.$$
\end{example}

\begin{example} Let $T(x)=\dis\frac{1}{1-x}+(q-1)x^r$, that is,
$t_i=1$ for all $i\geq 1$ except for $i=r$ and $t_r=q$. Then Theorem
\ref{theo 2.2} gives
\begin{eqnarray*}
\sum_{P\in\mathfrak{P}_{n,k}^{m,j}}q^{\alpha_r(P)}
&=&\frac{1}{m+1}\binom{m+j}{m}\sum_{i=0}^m\binom{m+j+1}{i}\binom{2m+j-(r+1)i}{m+j-i}(q-1)^i,
\end{eqnarray*}
which implies that the number of $k$-paths $P$ of length $n$ with
exactly $p$ $u$-segments of length $r$ $(namely\ \alpha_r(P)=p)$ is
given by
$\frac{1}{m+1}\binom{m+j}{m}\sum_{i=0}^m(-1)^{i-p}\binom{m+j+1}{i}\binom{2m+j-(r+1)i}{m+j-i}\binom{i}{p}$.
\end{example}

\begin{example} If $T(x)=1+\dis\frac{qx^r}{1-x^r}=\dis\frac{1+(q-1)x^r}{1-x^r}$, then
$t_{i}=q$ if $i\equiv 0\ mod\ r$ and $0$ otherwise. Thus, Theorem
\ref{theo 2.2} leads to  $$\begin{array}{ll}
\sum\limits_{P\in\mathfrak{P}_{2rm+kj,k}^{rm,j}}\prod\limits_{i\geq
1,\ i\equiv 0\ mod\ r}{q^{\alpha_{i}(P)}}
&=\frac{1}{rm+1}\binom{rm+j}{rm}\sum\limits_{i=1}^{m}\binom{rm+j+1}{i}\binom{m-1}{m-i}q^i\\
&=\frac{1}{rm+1}\binom{rm+j}{rm}\sum\limits_{i=0}^{m}\binom{rm+j+1}{i}\binom{(r+1)m+j-i}{m-i}(q-1)^i.
\end{array}$$
which produces the following results. The number of $k$-paths in
$\mathfrak{P}_{2rm+kj,k}^{rm,j}$ such that the length of any
$u$-segment is a multiple of $r$ $(i.e., the\ case\ $q=1$)$ is given
by $\frac{1}{rm+1}\binom{rm+j}{rm}\binom{(r+1)m+j}{m}$ (by
Vandermonde convolution). More precisely, the number of $k$-paths in
$\mathfrak{P}_{2rm+kj,k}^{rm,j}$ with exactly $i$ $u$-segments such
that the length of any $u$-segment is a multiple of $k$ is given by
$$\begin{array}{l}
\frac{1}{rm+1}\binom{rm+j}{rm}\binom{rm+j+1}{i}\binom{m-1}{m-i}=
\frac{1}{rm+1}\binom{rm+j}{rm}\sum\limits_{p=0}^{m}(-1)^{p-i}\binom{rm+j+1}{p}\binom{(r+1)m+j-p}{m-p}\binom{p}{i}.
\end{array}$$
\end{example}

\begin{example}\label{ex11}
Let $T(x)$ be the generating function $f^r(x)$, where $f(x)$ is the
generating function for complete $p$-ary plane trees (see, for
instance, \cite{chen99,klarner} and \cite[pp.112-113]{goulden}),
which satisfies the relation $f(x)=1+xf^p(x)$. By the Lagrange
inversion formula (\ref{eqn 1.2}), we can deduce
$t_i=\frac{r}{pi+r}\binom{pi+r}{i}$. Then Theorem \ref{theo 2.2}
leads to
$$\begin{array}{l}
\sum\limits_{P\in\mathfrak{P}_{n,k}^{m,j}}\prod\limits_{i\geq
1}\left\{\frac{r}{pi+r}\binom{pi+r}{i}\right\}^{\alpha_i(P)}
=\frac{1}{(m+1)}\frac{(m+j+1)r}{(m+j+1)r+mp}\binom{m+j}{m}\binom{r(m+j+1)+mp}{m}.
\end{array}$$
\end{example}

\subsection{A combinatorial proof of Corollary \ref{coro 2.3}}

Let $\mathfrak{\hat{P}}_{n,k,r}^{m,j}$ be the set of lattice paths
$P^*=Pd$ such that there is one colored down-step in $P^*$, where
$P\in \mathfrak{P}_{n,k,r}^{m,j}$. To give a bijective proof of
Corollary \ref{coro 2.3}, we need the following lemma.

\begin{lemma}\label{lemma 7.1}
There exists a bijection $\phi$ between the sets
$\mathfrak{\hat{P}}_{n,k,r}^{m,j}$ and $\mathfrak{L}_{n,k,r}^{m,j}$
such that $P^*\in \mathfrak{\hat{P}}_{n,k,r}^{m,j}$ has the same
type of $u$-segments as $\phi(P^*)\in \mathfrak{L}_{n,k,r}^{m,j}$.
\end{lemma}
\begin{proof}
Any $P^*\in \mathfrak{\hat{P}}_{n,k,r}^{m,j}$ can be uniquely
partitioned into $P^*=P_1dQ_1$, where $P_1, Q_1$ are lattice paths
and $d$ is the colored down-step. Define $\phi(P^*)=Q_1P_1$, then it
is easy to verify that $\phi(P^*)\in \mathfrak{L}_{n,k,r}^{m,j}$ and
note that the length of any $u$-segment in $\phi(P^*)$ is the same
as in $P^*$.

Conversely, for any $L\in \mathfrak{L}_{n,k,r}^{m,j}$, we can find
the leftmost point which has the lowest ordinate, then $L$ can be
uniquely partitioned into two parts in this sense, namely,
$L=L_1L_2$. Define $\phi^{-1}(L)=L_2dL_1$, where the $d$ is the
colored down-step, then it is easily to verify that $\phi^{-1}(L)\in
\mathfrak{\hat{P}}_{n,k,r}^{m,j}$ which has the same type of
$u$-segments with $L$.

Hence $\phi$ is indeed a bijection between the sets
$\mathfrak{\hat{P}}_{n,k,r}^{m,j}$ and $\mathfrak{L}_{n,k,r}^{m,j}$,
which preserves the type of $u$-segments not changed.
\end{proof}

An {\em ordered partition} $B_1,B_2,\cdots,B_r$ of
$[m]=\{1,2,\cdots,m\}$ into $r$ blocks is a partition of $[m]$ such
that the $r$ blocks as well as the elements of each block are
ordered.

Now we can give a bijective proof of Corollary \ref{coro 2.3}.
\begin{proof}
For any ordered partition $B_1,B_2,\cdots,B_r$ of $[m]$ into $r$
blocks, regard each block $B_i$ as a labeled $u$-segment $U_i$ for
$1\leq i\leq r$. For $m$ down-steps and $j$ horizontal-steps, there
are $\binom{m+j}{j}$ ways to obtain $(d,h)$-words of length $m+j$ on
$\{d,h\}$ with $m$ $d$'s and $j$ $h$'s. Then we can insert the
labeled $u$-segments $U_1,U_2,\cdots, U_r$ orderly into the $m+j+1$
positions (repetition is not allowed) of any $(d,h)$-word of length
$m+j$, which can produce $\binom{m+j+1}{r}$ labeled lattice paths in
$\mathfrak{L}_{n,k,r}^{m,j}$. Note that $r!{\bf
B}_{m,r}\big(1!t_1,2!t_2,\cdots\big)$ is just the generating
function for ordered partitions $B_1,B_2,\cdots,B_r$ of $[m]$ into
$r$ blocks such that each block $B_p$ is weighted by $t_i$ with
$i=|B_p|$ for $1\leq p\leq r$. So
${\binom{m+j}{m}}{\binom{m+j+1}{r}}r!{\bf
B}_{m,r}\big(1!t_1,2!t_2,\cdots\big)$ is the generating function for
the labeled lattice paths in $\mathfrak{L}_{n,k,r}^{m,j}$ such that
each $u$-segment of length $i$ is weighted by $t_i$.

However, by Lemma \ref{lemma 7.1}, any $k$-path $P\in
\mathfrak{P}_{n,k,r}^{m,j}$ can lead to $m!$ labeled $k$-paths, and
$P^*=Pd \in \mathfrak{\hat{P}}_{n,k,r}^{m,j}$ can generate $m+1$
lattice paths in $\mathfrak{L}_{n,k,r}^{m,j}$ and vice versa. Hence
$\frac{1}{m+1}{\binom{m+j}{m}}{\binom{m+j+1}{r}}\frac{r!}{m!}{\bf
B}_{m,r}\big(1!t_1,2!t_2,\cdots\big)$ is the generating function of
$k$-paths in $\mathfrak{P}_{n,k,r}^{m,j}$ such that each $u$-segment
of length $i$ is weighted by $t_i$, which makes the proof complete.
\end{proof}

\section{"internal $u$-segments" statistics in $k$-paths}

An {\em internal $u$-segment} of a $k$-path $P$ is a $u$-segment
between two steps such as $dd$, $hh$, $hd$, $dh$, i.e., all
$u$-segments except for the first one are internal $u$-segments.
Define $\beta_r(P)$ to be the number internal $u$-segments of length
$r$ in a $k$-path $P$. We start this section by studying the
ordinary generating functions for the number of $k$-paths of length
$n$ according to the statistics $\beta_1,\beta_2,\ldots$, that is,
\begin{eqnarray*}
F(x,z;{\bf t})=F(x,z;t_1,t_2,\ldots)=\sum_{m,j\geq
0}x^mz^j\sum_{P\in
\mathfrak{P}_{n,k}^{m,j}}\prod_{i\geq1}t_i^{\beta_i(P)},
\end{eqnarray*}
which can be represented as follows in terms of the generating
function $P(x,z;{\bf t})$.
\begin{proposition}\label{pro2}
The ordinary generating function $F(x,z;{\bf t})$ is given by
\begin{eqnarray*}
\frac{1+zP(x,z;{\bf t})}{1-xP(x,z;{\bf t})}=1+zP(x,z;{\bf
t})+\sum_{j\geq1}x^jP^j(x,z;{\bf
t})+z\sum_{j\geq1}x^jP^{j+1}(x,z;{\bf t}).
\end{eqnarray*}
\end{proposition}
\begin{proof}
An equation for $F(x,z;{\bf t})$ is obtained from the decomposition
of a $k$-path: either
$$P=hP',\ P=u^jdP^{(j)}dP^{(j-1)}\cdots dP^{(2)}dP^{(1)},\mbox{ or }
P=u^jhP^{(j+1)}dP^{(j)}\cdots dP^{(2)}dP^{(1)}$$ for some $j\geq 1$,
where $P', P^{(1)},\cdots,P^{(j+1)}$ are $k$-paths. Then $F(x,z;{\bf
t})$ satisfies the equation $F(x,z;{\bf t})=1+zP(x,z;{\bf
t})+\sum_{j\geq1}x^jP^j(x,z;{\bf
t})+z\sum_{j\geq1}x^jP^{j+1}(x,z;{\bf t})$, as required.
\end{proof}

\begin{theorem}\label{theo 3.2}
For any integers $k,j\geq 0$, $n,m\geq 1$,
\begin{eqnarray*}
\lefteqn{\sum_{P\in
\mathfrak{P}_{n,k}^{m,j}}\prod_{i\geq1}t_i^{\beta_i(P)}=
\sum_{p=0}^{m}\left\{\frac{m-p}{m+j}\binom{m+j}{j}+\binom{m+j}{j-1}\frac{m-p+1}{m+j}\right\}{\bf
P}_{p}^{(m+j)}\big(1!t_1,2!t_2,\cdots\big)}\\
&=&\binom{m+j}{j}\sum_{p=0}^{m}\frac{(m+1)(m-p)+j(m-p+1)}{(m+1)(m+j)}\sum_{r=0}^{p}{\binom{m+j}{r}}\frac{r!}{p!}{\bf
B}_{p,r}\big(1!t_1,2!t_2,\cdots\big).
\end{eqnarray*}
\end{theorem}
\begin{proof}
Using (\ref{eqn 1.2}) and (\ref{eqn 1.3}), we obtain
\begin{eqnarray*}
\lefteqn{\sum_{P\in
\mathfrak{P}_{n,k}^{m,j}}\prod_{i\geq1}t_i^{\beta_i(P)}=
[x^{m}z^j]\frac{1+zP(x,z;{\bf t})}{1-xP(x,z;{\bf
t})}=[x^{m+j}z^j]\frac{1+zy^*}{1-y^*}}\\
&=&\frac{1}{m+j}[x^{m+j-1}z^j]\left\{\frac{1+zx}{1-x}\right\}'\big((1+zx)T(x)\big)^{m+j}\\
&=&\frac{1}{m+j}\binom{m+j}{j}[x^{m-1}]\frac{T(x)^{m+j}}{(1-x)^2}
+\frac{1}{m+j}\binom{m+j}{j-1}[x^{m}]\frac{T(x)^{m+j}}{(1-x)^2}\\
&=&\binom{m+j}{j}\sum_{p=0}^{m}\frac{m-p}{m+j}{\bf
P}_{p}^{(m+j)}\big(1!t_1,2!t_2,\cdots\big)
+\binom{m+j}{j-1}\sum_{p=0}^{m}\frac{m-r+1}{m+j}{\bf
P}_{p}^{(m+j)}\big(1!t_1,2!t_2,\cdots\big)\\
&=&\sum_{p=0}^{m}\left\{\frac{m-p}{m+j}\binom{m+j}{j}+\binom{m+j}{j-1}\frac{m-p+1}{m+j}\right\}{\bf
P}_{p}^{(m+j)}\big(1!t_1,2!t_2,\cdots\big)\\
&=&\binom{m+j}{j}\sum_{p=0}^{m}\frac{(m+1)(m-p)+j(m-p+1)}{(m+1)(m+j)}\sum_{r=0}^{p}{\binom{m+j}{r}}\frac{r!}{p!}{\bf
B}_{p,r}\big(1!t_1,2!t_2,\cdots\big),
\end{eqnarray*}
as claimed.
\end{proof}

\section{"$u$-segments" and "Internal $u$-segments" Statistics in $k$-paths without a horizontal-step on the $x$-axis}
\subsection{$u$-segments statistics} We start this subsection by studying the generating function for the
number of $k$-paths in $\mathfrak{Q}_{n,k}^{m,j}$ according to the
statistics $\alpha_1,\alpha_2,\ldots$, that is,
\begin{eqnarray*}
Q(x,z;{\bf
t})=Q(x,z;t_1,t_2,\ldots)=\sum_{m,j\geq0}x^mz^j\sum_{Q\in\mathfrak{Q}_{n,k}^{m,j},
}\prod_{i\geq1}t_i^{\alpha_i(Q)}.
\end{eqnarray*}

\begin{proposition}\label{pro4.1}
The ordinary generating function $Q(x,z;{\bf t})$ is given by
\begin{eqnarray}\label{eqn 4.1}
Q(x,z;{\bf t})=\frac{P(x,z;{\bf t})}{1+zP(x,z;{\bf t})}=T(y).
\end{eqnarray}
\end{proposition}
\begin{proof}
Note that $Q(x,z;{\bf t})$ can be written as $Q(x,z;{\bf
t})=1+\sum_{j\geq1}Q_p(x,z;{\bf t})$, where $Q_p(x,z;{\bf t})$ is
the generating function for the number of $k$-paths starting with
$p$ up-steps and without a horizontal-step on the $x$-axis according
to the statistics $\alpha_1,\alpha_2,\ldots$. An equation for
$Q_p(x,z;{\bf t})$ is obtained from the first return decomposition
of a $k$-path starting with a $u$-segment of length $p$: either
$P=u^pdP^{(p-1)}d P^{(p-2)} d\ldots P^{(1)}d P^{*}$ or $P=
u^phP^{(p)} d P^{(p-1)} d\ldots P^{(1)}dP^{*}$, where
$P^{(1)},\ldots,P^{(p)}$ are $k$-paths and $P^{*}$ is a $k$-path
without a horizontal-step on the $x$-axis. Thus $Q_p(x,z;{\bf
t})=t_{p}x^{p}P^{p-1}(x,z;{\bf t})Q(x,z;{\bf
t})+zt_{p}x^{p}P^{p}(x,z;{\bf t})Q(x,z;{\bf t})$ and $Q(x,z;{\bf
t})$ satisfies the equation $Q(x,z;{\bf t})=1+Q(x,z;{\bf
t})\big\{\sum_{p\geq1}t_px^{p}P^{p-1}(x,z;{\bf
t})+z\sum_{p\geq1}t_px^{p}P^{p}(x,z;{\bf t})\big\}$. Hence, by
Proposition \ref{pro2.1}, we obtain the desired result.
\end{proof}

\begin{theorem}\label{theo 4.2}
For any integers $k,j\geq 0$, $n,m\geq 1$,
\begin{eqnarray*}
\lefteqn{\sum_{Q\in
\mathfrak{Q}_{n,k}^{m,j}}\prod_{i\geq1}t_i^{\alpha_i(Q)}
=\frac{m}{(m+j+1)(m+j)}\binom{m+j}{j}\frac{{\bf{P}}_m^{(m+j+1)}(1!t_1,2!t_2,\cdots)}{m!}}\\
&=&\frac{m}{(m+j+1)(m+j)}\binom{m+j}{j}\sum_{r=0}^m{\binom{m+j+1}{r}}\frac{r!}{m!}{\bf
B}_{m,r}\big(1!t_1,2!t_2,\cdots\big).
\end{eqnarray*}
\end{theorem}
\begin{proof}
Using (\ref{eqn 1.2}) and (\ref{eqn 1.3}), we obtain
\begin{eqnarray*}
\lefteqn{\sum_{Q\in
\mathfrak{Q}_{n,k}^{m,j}}\prod_{i\geq1}t_i^{\alpha_i(Q)}=[x^mz^j]Q(x,z;{\bf
t})=[x^{m+j}z^j]T(y^*)}\\
&=&\frac{1}{m+j}[x^{m+j-1}z^j]\left\{T(x)\right\}'\left\{(1+zx)T(x)\right\}^{m+j}\\
&=&\frac{1}{m+j}\binom{m+j}{j}[x^{m-1}]\left\{T(x)\right\}'\left\{T(x)\right\}^{m+j}\\
&=&\frac{1}{m+j}\binom{m+j}{j}\frac{1}{m+j+1}[x^{m-1}]\left\{T(x)^{m+j+1}\right\}'\\
&=&\frac{m}{(m+j+1)(m+j)}\binom{m+j}{j}[x^{m}]T(x)^{m+j+1}\\
&=&\frac{m}{(m+j+1)(m+j)}\binom{m+j}{j}\frac{{\bf{P}}_m^{(m+j+1)}(1!t_1,2!t_2,\cdots)}{m!}\\
&=&\frac{m}{(m+j+1)(m+j)}\binom{m+j}{j}\sum_{r=0}^m{\binom{m+j+1}{r}}\frac{r!}{m!}{\bf
B}_{m,r}\big(1!t_1,2!t_2,\cdots\big),
\end{eqnarray*}
which completes the proof.
\end{proof}

Replacing $t_i$ by $qt_i$ in Theorem \ref{theo 4.2} and using
(\ref{eqn bella}) and (\ref{eqn bellb}), we have
\begin{corollary}\label{coro 4.3}
For any integers $n,m,r\geq 1$ and $k,j\geq 0$, there holds
\begin{eqnarray*}\label{eqn 4.2}
\sum_{Q\in\mathfrak{Q}_{n,k,r}^{m,j}}\prod_{i\geq1}t_i^{\alpha_i(Q)}
&=&\frac{m}{(m+j+1)(m+j)}\binom{m+j}{j}{\binom{m+j+1}{r}}\frac{r!}{m!}{\bf
B}_{m,r}\big(1!t_1,2!t_2,\cdots\big),
\end{eqnarray*}
which implies that the number of $k$-paths with no horizontal-step
on the $x$-axis and with $u$-segments of type $1^{r_1}2^{r_2}\cdots
m^{r_m}$ is
$\frac{m}{(m+j+1)(m+j)}{\binom{m+j}{m}}{\binom{m+j+1}{r}}\binom{r}{r_1,r_2,\cdots,r_m}$.
\end{corollary}

\begin{remark} Note that from Theorem \ref{theo 2.2} and \ref{theo 4.2}, the ratio
of
$\sum_{P\in\mathfrak{P}_{n,k}^{m,j}}\prod_{i\geq1}t_i^{\alpha_i(P)}$
to
$\sum_{Q\in\mathfrak{Q}_{n,k}^{m,j}}\prod_{i\geq1}t_i^{\alpha_i(Q)}$
is $\frac{(m+1)m}{(m+j+1)(m+j)}$. For the sake of conciseness, we
omit many examples as done in Section 2. However, we may ask whether
there is a combinatorial interpretation for this relation.
\end{remark}

\subsection{Internal $u$-segments statistics}
In this subsection, we study the generating function for the number
of $k$-paths in $\mathfrak{Q}_{n,k}^{m,j}$ according to the
statistics $\beta_1,\beta_2,\ldots$, that is,
\begin{eqnarray*}
H(x,z;{\bf
t})=H(x,z;t_1,t_2,\ldots)=\sum_{m,j\geq0}x^mz^j\sum_{Q\in\mathfrak{Q}_{n,k}^{m,j},
}\prod_{i\geq1}t_i^{\beta_i(Q)}.
\end{eqnarray*}
\begin{proposition}\label{pro5.1}
The ordinary generating function $H(x,z;{\bf t})$ is given by
\begin{eqnarray}\label{eqn 5.1}
H(x,z;{\bf t})=\frac{1-xP(x,z;{\bf t})}{1-x-xP(x,z;{\bf
t})-zxP(x,z;{\bf t})}=\frac{1-y}{1-y-x-zy}.
\end{eqnarray}
\end{proposition}
\begin{proof}
Note that $H(x,z;{\bf t})$ can be written as $H(x,z;{\bf
t})=1+\sum_{p\geq1}H_p(x,z;{\bf t})$, where $H_p(x,z;{\bf t})$ is
the generating function for the number of $k$-paths starting with
$p$ up-steps and without a horizontal-step on the $x$-axis according
to the statistics $\beta_1,\beta_2,\ldots$. An equation for
$H_p(x,z;{\bf t})$ is obtained from the first return decomposition
of a $k$-path starting with a $u$-segment of length $p$: either
$P=u^pdP^{(p-1)}d P^{(p-2)} d\ldots P^{(1)}d P^{*}$ or $P=
u^phP^{(p)} d P^{(p-1)} d\ldots P^{1}dP^{*},$ where
$P^{(1)},\ldots,P^{(p)}$ are $k$-paths and $P^{*}$ is a $k$-path
with no horizontal-step on the $x$-axis. Thus $H_p(x,z;{\bf
t})=x^{p}P^{p-1}(x,z;{\bf t})H(x,z;{\bf t})+zx^{p}P^{p}(x,z;{\bf
t})H(x,z;{\bf t})$. Hence, $H(x,z;{\bf t})$ satisfies the equation
$H(x,z;{\bf t})=1+\big\{\sum_{p\geq1}x^{p}P^{p-1}(x,z;{\bf
t})+z\sum_{p\geq1}x^{p}P^{p}(x,z;{\bf t})\big\}H(x,z;{\bf t})$, a
simplification reduces this to the required expression.
\end{proof}

\begin{theorem}\label{theo 5.2}
For any integers $n,m,k,j\geq 0$, $m+j\geq 1$,
\begin{eqnarray*}
\lefteqn{\sum_{Q\in
\mathfrak{Q}_{n,k}^{m,j}}\prod_{i\geq1}t_i^{\beta_i(Q)}=\sum_{i=0}^m
\frac{i}{m+j-i}\binom{m+j-1}{j}\sum_{r=0}^{m-i-1}\binom{r+i}{r}
\frac{\textbf{P}_{m-i-r-1}^{(m+j-i)}(1!t_1,2!t_2,\cdots)}{(m-i-r-1)!}}\\
&\hskip2.5cm&+\sum_{i=0}^m
\frac{i}{m+j-i}\binom{m+j-1}{j-1}\sum_{r=0}^{m-i}\binom{r+i}{r}
\frac{\textbf{P}_{m-i-r}^{(m+j-i)}(1!t_1,2!t_2,\cdots)}{(m-i-r)!}.
\end{eqnarray*}
\end{theorem}
\begin{proof}
Using (\ref{eqn 1.2}) and (\ref{eqn 1.3}), we obtain
\begin{eqnarray*}
\lefteqn{\sum_{Q\in
\mathfrak{Q}_{n,k}^{m,j}}\prod_{i\geq1}t_i^{\beta_i(Q)}=[x^mz^j]Q(x,z;{\bf
t})=[x^{m+j}z^j]\frac{1-y^*}{1-y^*-x(1+zy^*)}}\\
&=&[x^{m+j}z^j]\sum_{i\geq
0}^{m+j}\left\{\frac{x(1+zy^*)}{1-y^*}\right\}^i
=\sum_{i\geq 0}^{m+j}[x^{m+j-i}z^j]\left\{\frac{1+zy^*}{1-y^*}\right\}^i\\
&=&\sum_{i\geq
0}^{m+j}\frac{i}{m+j-i}[x^{m+j-i-1}z^j]\left\{\frac{1+zx}{1-x}\right\}^{i-1}\left\{\frac{1+zx}{1-x}\right\}'
\left\{(1+zx)T(x)\right\}^{m+j-i}\\
&=&\sum_{i\geq
0}^{m+j}\frac{i}{m+j-i}[x^{m+j-i-1}z^j]\left\{\frac{(1+z)(1+zx)^{m+j-1}}{(1-x)^{i+1}}\right\}
\left\{T(x)\right\}^{m+j-i}\\
&=&\sum_{i\geq
0}^{m+j}\frac{i}{m+j-i}\left\{\binom{m+j-1}{j}[x^{m-i-1}]\frac{T(x)^{m+j-i}}{(1-x)^{i+1}}+
\binom{m+j-1}{j-1}[x^{m-i}]\frac{T(x)^{m+j-i}}{(1-x)^{i+1}}\right\}\\
&=&\sum_{i=0}^m
\frac{i}{m+j-i}\binom{m+j-1}{j}\sum_{r=0}^{m-i-1}\binom{r+i}{r}
\frac{\textbf{P}_{m-i-r-1}^{(m+j-i)}(1!t_1,2!t_2,\cdots)}{(m-i-r-1)!}\\
&\hskip0.5cm&+\sum_{i=0}^m
\frac{i}{m+j-i}\binom{m+j-1}{j-1}\sum_{r=0}^{m-i}\binom{r+i}{r}
\frac{\textbf{P}_{m-i-r}^{(m+j-i)}(1!t_1,2!t_2,\cdots)}{(m-i-r)!},
\end{eqnarray*}
which completes the proof.
\end{proof}

\section{ Statistics $(u,h)$-segments and  $u$-segments in $k$-paths}
A {\em{$(u,h)$-segment}} in a $k$-path is a maximum segment composed
of up-steps and horizontal-steps. An {\em{internal $(u,h)$-segment}}
in a $k$-path is a {$(u,h)$-segment} between two down steps. Let
$\mathfrak{\tilde{P}}_{n,k}^{m,j,\ell}$ denote the subset of
$\mathfrak{P}_{n,k}^{m,j}$ such that $\rm (i)$ each internal
$(u,h)$-segment has length equal to a multiple of $k$; $\rm (ii)$
the first $(u,h)$-segment has length $\equiv \ell\ (mod\ k)$ for
$0\leq \ell\leq k-1$. We note that the case $j=0$ is studied in
\cite{mansun}.

\begin{theorem}\label{theo 6.1}
The number $\mathcal{\tilde{P}}_{r,k,\ell}$ of $k$-paths of length
$n=kr+2\ell$ satisfying the conditions $\rm (i)$ and $\rm (ii)$ is
\begin{eqnarray*}
\mathcal{\tilde{P}}_{r,k,\ell}
&=&\sum_{p=0}^{[r/2]}\frac{\ell+1}{kp+\ell+1}\binom{(k+1)p+\ell}{p}\binom{r+2(k-1)p+2\ell}{r-2p}.
\end{eqnarray*}
\end{theorem}
\begin{proof}
Note that for any path $P\in \mathfrak{\tilde{P}}_{n,k}^{m,j,\ell}$,
by deleting all the $j$ horizontal-steps, we obtain a Dyck path in
$\mathfrak{\tilde{P}}_{2m,k}^{m,0,\ell}$. Conversely, any Dyck path
in $\mathfrak{\tilde{P}}_{2m,k}^{m,0,\ell}$ goes through $2m+1$
integer points, if we insert the $j$ horizontal-steps into any
integer point (repetitions are allowed), then we get
$\binom{2m+j}{j}$ $k$-paths in
$\mathfrak{\tilde{P}}_{n,k}^{m,j,\ell}$. However, the set
$\mathfrak{\tilde{P}}_{2m,k}^{m,0,\ell}$ is counted by
$\frac{\ell+1}{m+1}\binom{m+p}{p}$, which has been proved in
\cite{mansun}. Hence we have
\begin{eqnarray*}
\mathcal{\tilde{P}}_{r,k,\ell}&=&\sum_{2m+kj=n,m=kp+\ell}{|\mathfrak{\tilde{P}}_{n,k}^{m,j,\ell}|}\\
&=&\sum_{2m+kj=n,m=kp+\ell}\frac{\ell+1}{m+1}\binom{m+p}{p}\binom{2m+j}{j}\\
&=&\sum_{p=0}^{[r/2]}\frac{\ell+1}{kp+\ell+1}\binom{(k+1)p+\ell}{p}\binom{r+2(k-1)p+2\ell}{r-2p},
\end{eqnarray*}
as required.
\end{proof}

Let $\mathfrak{\bar{P}}_{n,k}^{m,j,\ell}$ denote the subset of
$\mathfrak{P}_{n,k}^{m,j}$ such that $\rm (iii)$ each internal
$u$-segment has length equal to a multiple of $k$; $\rm (iv)$ the
first $u$-segment has length $\equiv \ell\ (mod\ k)$ for $0\leq
\ell\leq k-1$.

\begin{theorem}\label{theo 6.2}
The number $\mathcal{\bar{P}}_{r,k,\ell}$ of $k$-paths of length
$n=kr+2\ell$ satisfying conditions $\rm (iii)$ and $\rm (iv)$ is
\begin{eqnarray*}
\mathcal{\bar{P}}_{r,k,\ell}
&=&\sum_{p=0}^{[r/2]}\frac{\ell+1}{kp+\ell+1}\binom{(k+1)p+\ell}{p}\binom{r+(k-1)p+\ell}{r-2p}.
\end{eqnarray*}
\end{theorem}
\begin{proof}
Note that for any path $P\in\mathfrak{\bar{P}}_{n,k}^{m,j,\ell}$, by
deleting all the $j$ horizontal-steps, we obtain a Dyck path in $
\mathfrak{\bar{P}}_{2m,k}^{m,0,\ell}$. Conversely,let $D$ be a Dyck
path from $\mathfrak{\bar{P}}^{m,0,\ell}_{2m,k}$, where $m=kp+\ell$
for some $p\geq 0$. It can be shown that $D$ has $m+p+1$ proper
integer points, where a proper integer point is a point where we may
insert a horizontal step without violating the properties (iii) and
(iv). By inserting $j$ horizontal steps into these points
(repetitions are allowed) we get $\binom{m+p+j}{j}$ $k$-paths in
$\mathfrak{\bar{P}}^{m,j,\ell}_{2m,k}$. Note that the set
$\mathfrak{\bar{P}}_{2m,k}^{m,0,\ell}=\mathfrak{\tilde{P}}_{2m,k}^{m,0,\ell}$
is counted by $\frac{\ell+1}{m+1}\binom{m+p}{p}$. Hence we have
\begin{eqnarray*}
\mathcal{\bar{P}}_{r,k,\ell}&=&\sum_{2m+kj=n,m=kp+\ell}{|\mathfrak{\bar{P}}_{n,k}^{m,j,\ell}|}\\
&=&\sum_{2m+kj=n,m=kp+\ell}\frac{\ell+1}{m+1}\binom{m+p}{p}\binom{m+p+j}{j}\\
&=&\sum_{p=0}^{[r/2]}\frac{\ell+1}{kp+\ell+1}\binom{(k+1)p+\ell}{p}\binom{r+(k-1)p+\ell}{r-2p},
\end{eqnarray*}
as required.
\end{proof}

It should be noted that $\sum_{p\geq
0}\frac{\ell+1}{kp+\ell+1}\binom{(k+1)p+\ell}{p}x^p=f(x)^{\ell+1}$,
where $f(x)$ is the generating function for $(k+1)$-ary plane trees,
and which satisfies the relation $f(x)=1+xf(x)^{k+1}$. Then it is
easy to prove that the generating functions for
$\mathcal{\tilde{P}}_{r,k,\ell}$ and for
$\mathcal{\bar{P}}_{r,k,\ell}$ are respectively

\begin{eqnarray}
\tilde{P}_{k,\ell}(x)&=&\sum_{r\geq
0}\mathcal{\tilde{P}}_{r,k,\ell}x^{r}=\frac{1}{(1-x)^{2\ell+1}}f(\frac{x^2}{(1-x)^{2k}})^{\ell+1}, \label{eqn 6.1}\\
\bar{P}_{k,\ell}(x)&=&\sum_{r\geq
0}\mathcal{\bar{P}}_{r,k,\ell}x^{r}=\frac{1}{(1-x)^{\ell+1}}f(\frac{x^2}{(1-x)^{k+1}})^{\ell+1}.
\label{eqn 6.2}
\end{eqnarray}
Replacing $x$ by $\frac{x}{1+x}$ in (\ref{eqn 6.1}) and (\ref{eqn
6.2}), one can deduce that
\begin{eqnarray}
\frac{1}{(1+x)^{2\ell+1}}\tilde{P}_{k,\ell}(\frac{x}{1+x})&=&f(x^2(1+x)^{2k-2})^{\ell+1}, \label{eqn 6.3}\\
\frac{1}{(1+x)^{\ell+1}}\bar{P}_{k,\ell}(\frac{x}{1+x})&=&f(x^2(1+x)^{k-1})^{\ell+1}.
\label{eqn 6.4}
\end{eqnarray}
Comparing the coefficient of $x^n$ in both sides of (\ref{eqn 6.3})
and (\ref{eqn 6.4}), one can deduce the following consequence:
\begin{corollary}\label{coro 6.0}
For any integers $n,k,\ell\geq 0$, there hold
\begin{eqnarray}
\sum_{p=0}^{n}(-1)^{n-p}\binom{n+2\ell}{p+2\ell}\mathcal{\tilde{P}}_{p,k,\ell}
&=&\sum_{p=0}^{[\frac{n}{2}]}\frac{\ell+1}{kp+\ell+1}\binom{(k+1)p+\ell}{p}\binom{2(k-1)p}{n-2p},\nonumber\\
\sum_{p=0}^{n}(-1)^{n-p}\binom{n+\ell}{p+\ell}\mathcal{\bar{P}}_{p,k,\ell}
&=&\sum_{p=0}^{[\frac{n}{2}]}\frac{\ell+1}{kp+\ell+1}\binom{(k+1)p+\ell}{p}\binom{(k-1)p}{n-2p}.
\label{eqn 6.5}
\end{eqnarray}
\end{corollary}
Using the generalized Lagrange inversion formula obtained in
\cite{corsani}, from (\ref{eqn 6.1}), we have
\begin{eqnarray*}
\lefteqn{\frac{\ell+1}{kn+\ell+1}\binom{(k+1)n+\ell}{n}=[x^n]f(x)^{\ell+1}}\\
&=&[w^n]\left\{(1-x)^{2\ell+1}\tilde{P}_{k,\ell}(x)\right\}_{w=\frac{x^2}{(1-x)^{2k}}}\\
&=&\frac{1}{2n}[t^{2n-1}](1-t)^{2kn}\frac{d}{dt}\left\{(1-t)^{2\ell+1}\tilde{P}_{k,\ell}(t)\right\}\\
&=&\frac{2\ell+1}{2n}[t^{2n-1}](1-t)^{2kn+2\ell}\tilde{P}_{k,\ell}(t)
+\frac{1}{2n}[t^{2n-1}](1-t)^{2kn+2\ell+1}\frac{d}{dt}\tilde{P}_{k,\ell}(t)\\
&=&\sum_{p=0}^{2n-1}(-1)^p\binom{2kn+2\ell}{2n-p-1}\frac{2\ell+1}{2n}\mathcal{\tilde{P}}_{p,k,\ell}
+\sum_{p=1}^{2n}(-1)^p\binom{2kn+2\ell+1}{2n-p}\frac{p}{2n}\mathcal{\tilde{P}}_{p,k,\ell}\\
&=&\sum_{p=0}^{2n}(-1)^p\frac{kp+2\ell+1}{2kn+2\ell+1}\binom{2kn+2\ell+1}{2n-p}\mathcal{\tilde{P}}_{p,k,\ell}.
\end{eqnarray*}
Similarly, from (\ref{eqn 6.2}), we have
\begin{eqnarray*}
\lefteqn{\frac{\ell+1}{kn+\ell+1}\binom{(k+1)n+\ell}{n}=[x^n]f(x)^{\ell+1}}\\
&=&[w^n]\left\{(1-x)^{\ell+1}\bar{P}_{k,\ell}(x)\right\}_{w=\frac{x^2}{(1-x)^{k+1}}}\\
&=&\frac{1}{2n}[t^{2n-1}](1-t)^{(k+1)n}\frac{d}{dt}\left\{(1-t)^{\ell+1}\bar{P}_{k,\ell}(t)\right\}\\
&=&\sum_{p=0}^{2n}(-1)^{p}\frac{p(k+1)+2(\ell+1)}{2n(k+1)+2(\ell+1)}\binom{n(k+1)+\ell+1}{2n-p}\mathcal{\bar{P}}_{p,k,\ell}.
\end{eqnarray*}
Hence we obtain the next corollary:
\begin{corollary}\label{coro 6.1}
For any integers $n,k,\ell\geq 0$,it holds that
\begin{eqnarray}
\hskip0.6cm\frac{\ell+1}{kn+\ell+1}\binom{(k+1)n+\ell}{n}
&=&\sum_{p=0}^{2n}(-1)^p\frac{kp+2\ell+1}{2kn+2\ell+1}\binom{2kn+2\ell+1}{2n-p}\mathcal{\tilde{P}}_{p,k,\ell},\label{eqn 6.6}\\
\frac{\ell+1}{kn+\ell+1}\binom{(k+1)n+\ell}{n}&=&\sum_{p=0}^{2n}(-1)^{p}\frac{p(k+1)+2(\ell+1)}{2n(k+1)+2(\ell+1)}\binom{n(k+1)+\ell+1}{2n-p}
\mathcal{\bar{P}}_{p,k,\ell}.\label{eqn 6.7}
\end{eqnarray}
\end{corollary}

We consider below several special cases, leading to several
interesting results.
\begin{example} If $k=1$ and $\ell=0$ in (\ref{eqn 6.2}), then
$f(x)=\frac{1-\sqrt{1-4x}}{2x}=C(x)$,  which is the generating
function for the Catalan numbers $C_n=\frac{1}{n+1}\binom{2n}{n}$.
Hence we have
\begin{eqnarray*}
\bar{P}_{1,0}(x)=\frac{1}{1-x}C(\frac{x^2}{(1-x)^2})=\frac{1-x-\sqrt{1-2x-3x^2}}{2x^2},
\end{eqnarray*}
which is the generating function $M(x)$ for the Motzkin numbers
$M_n$. Then Theorem \ref{theo 6.2} together with (\ref{eqn 6.5}) and
(\ref{eqn 6.6}) generates the well-known identities (see
\cite{aigner,bernhart})
\begin{eqnarray*}
\sum_{p=0}^{[n/2]}\binom{n}{2p}C_{p}&=& M_n, \\
\sum_{p=0}^{n}(-1)^{n-p}\binom{n}{p}M_p&=&\left\{
\begin{array}{ll}
C_r & { \rm if}\ n=2r, \\
0   &  { \rm otherwise.}
\end{array}\right.
\end{eqnarray*}
\end{example}

\begin{example}\label{ex 6.1} If $k=2$ and $\ell=0$ in (\ref{eqn 6.2}), then
$f(x)=1+xf(x)^3=\sum_{n\geq 0}\frac{1}{2n+1}\binom{3n}{n}x^n$, which
is the generating function for complete $3$-ary plane trees. So we
have
\begin{eqnarray*}
\bar{P}_{2,0}(x)=\frac{1}{1-x}f(\frac{x^2}{(1-x)^3})
=\frac{1}{1-x}+\frac{x^2}{(1-x)^4}f(\frac{x^2}{(1-x)^3})^3=\frac{1}{1-x}(1+x^2\bar{P}_{2,0}(x)^3).
\end{eqnarray*}
If we let $y=x\bar{P}_{2,0}(x)$, it follows that $y=x(1+y)^{-1}$.
From this,
\begin{eqnarray*}
\bar{P}_{2,0}(x)=C(x)=\frac{1-\sqrt{1-4x}}{2x}.
\end{eqnarray*}
Similarly, when $k=2$ and $\ell=1$ in (\ref{eqn 6.2}), we have
\begin{eqnarray*}
\bar{P}_{2,1}(x)=\bar{P}_{2,0}(x)^2=C(x)^2=\frac{1-2x-\sqrt{1-4x}}{2x^2}=\frac{C(x)-1}{x}.
\end{eqnarray*}
\end{example}

Hence we obtain the following statement:

\begin{corollary}\label{coro 6.2}
The number of $2$-paths (i.e. Schr$\ddot{o}$der paths) of length
$2n$ such that all $u$-segments have even length is the Catalan
number $C_n$ for $n\geq 0$ and the number of $2$-paths of length
$2n+2$ such that all internal $u$-segments have even length and the
first $u$-segment has odd length is the Catalan number $C_{n+1}$ for
$n\geq 0$.
\end{corollary}

Here is a simple bijective proof. For any Schr\"oder path $S$ of
length $2n$ such that all $u$-segments have even length, replace
each $h$ step by $ud$ steps, then we get a Dyck path of length $2n$.
On the other hand, any Dyck path $D$ of length $2n$ can be
decomposed uniquely into $D=u^{i_1}d^{j_1}u^{i_2}d^{j_2}\cdots
u^{i_k}d^{j_k}$, where $i's,j's\geq 1$. Now replace a sub-path
$u^{i_l}d^{j_{l}}$ by $u^{i_l-1}hd^{j_{l}-1}$ if $i_l$ is odd, and
do nothing is $i_l$ is even. Then we get a desired Schr\"oder path
$S$. A similar argument proves the second claim in Corollary
\ref{coro 6.2}.

Theorem \ref{theo 6.2} and Example \ref{ex 6.1} give rise to two new
expressions for the Catalan numbers,

\begin{corollary}
For any integer $n\geq 0$,
\begin{eqnarray*}
\sum_{p=0}^{[n/2]}\frac{1}{2p+1}\binom{3p}{p}\binom{n+p}{3p}&=&\frac{1}{n+1}\binom{2n}{n},\hskip0.5cm (n\geq 0),\\
\sum_{p=0}^{[(n-1)/2]}\frac{1}{2p+1}\binom{3p+1}{p+1}\binom{n+p}{3p+1}&=&\frac{1}{n+1}\binom{2n}{n},\hskip0.5cm
(n\geq 1).
\end{eqnarray*}
\end{corollary}
Example \ref{ex 6.1} together with (\ref{eqn 6.5}) and (\ref{eqn
6.6}) leads to several new identities involving Catalan numbers
\begin{corollary}
For any integer $n\geq 0$,
\begin{eqnarray}
\sum_{p=0}^{n}(-1)^{n-p}\binom{n}{p}C_p&=&\sum_{p=0}^{[\frac{n}{2}]}\frac{1}{2p+1}\binom{3p}{p}\binom{p}{n-2p}, \label{eqn 6.8}\\
\sum_{p=0}^{n}(-1)^{n-p}\binom{n+1}{p+1}C_{p+1}&=&\sum_{p=0}^{[\frac{n}{2}]}\frac{1}{2p+1}\binom{3p+1}{p+1}\binom{p}{n-2p},\nonumber\\
\sum_{p=0}^{2n}(-1)^{p}\frac{3p+2}{2n+2p+2}\binom{3n}{n+p}C_p&=&\frac{1}{2n+1}\binom{3n}{n},\nonumber\\
\sum_{p=0}^{2n}(-1)^{p}\frac{3p+4}{2n+2p+4}\binom{3n+1}{n+p+1}C_{p+1}&=&\frac{1}{2n+1}\binom{3n+1}{n+1}.\nonumber
\end{eqnarray}
\end{corollary}
\begin{remark} $\rm (I)$ In fact, the counting formula in Theorem \ref{theo
6.1} and (\ref{eqn 6.6}), and the counting formula in Theorem
\ref{theo 6.2} and (\ref{eqn 6.7}) form two left-inversion relations
which have the general formats obtained implicitly by Corsani,
Merlini and Sprugnoli \cite{corsani}, namely
\begin{eqnarray*}
A_n=\sum_{p=0}^{[\frac{n}{2}]}\binom{n+2(k-1)p+2\ell}{n-2p}B_p
\Rightarrow
B_n=\sum_{p=0}^{2n}(-1)^{p}\frac{kp+2\ell+1}{2kn+2\ell+1}\binom{2kn+2\ell+1}{2n-p}A_p,\\
A_n=\sum_{p=0}^{[\frac{n}{2}]}\binom{n+(k-1)p+\ell}{n-2p}B_p
\Rightarrow
B_n=\sum_{p=0}^{2n}(-1)^{p}\frac{p(k+1)+2(\ell+1)}{2n(k+1)+2(\ell+1)}\binom{n(k+1)+\ell+1}{2n-p}A_p.
\end{eqnarray*}
$\rm (II)$ It should be noticed that the left side of (\ref{eqn
6.8}) is the Riordan numbers obtained by Bernhart \cite{bernhart},
so the right side of (\ref{eqn 6.8}) gives a new expression for the
Riordan numbers.
\end{remark}\vskip0.5cm

\subsection*{Acknowledgements} The authors are deeply grateful to the
two anonymous referees for valuable suggestions on an earlier
version of this paper which makes it more readable. Thanks also to
Simone Severini for helpful comments. The second author is supported
by (NSFC10726021) The National Science Foundation of China.


\end{document}